\documentclass[12pt]{amsart}

\usepackage{latexsym}
\usepackage{amsfonts}
\usepackage{amssymb}
\usepackage{amsmath}
\usepackage{mathrsfs}
\input xypic
\usepackage{hyperref}

\newcommand{\be}{\begin{equation}}
\newcommand{\ee}{\end{equation}}
\newcommand{\bea}{\begin{eqnarray}}
\newcommand{\eea}{\end{eqnarray}}
\newcommand{\bean}{\begin{eqnarray*}}
\newcommand{\eean}{\end{eqnarray*}}
\newcommand{\brray}{\begin{array}}
\newcommand{\erray}{\end{array}}

\newtheorem{dfn}{Definition}[section]
\newtheorem{thm}[dfn]{Theorem}
\newtheorem{lmma}[dfn]{Lemma}
\newtheorem{ppsn}[dfn]{Proposition}
\newtheorem{crlre}[dfn]{Corollary}
\newtheorem{xmpl}[dfn]{Example}
\newtheorem{rmrk}[dfn]{Remark}
\newcommand{\bdfn}{\begin{dfn}\rm}
\newcommand{\bthm}{\begin{thm}}
\newcommand{\blmma}{\begin{lmma}}
\newcommand{\bppsn}{\begin{ppsn}}
\newcommand{\bcrlre}{\begin{crlre}}
\newcommand{\bxmpl}{\begin{xmpl}}
\newcommand{\brmrk}{\begin{rmrk}\rm}

\newcommand{\edfn}{\end{dfn}}
\newcommand{\ethm}{\end{thm}}
\newcommand{\elmma}{\end{lmma}}
\newcommand{\eppsn}{\end{ppsn}}
\newcommand{\ecrlre}{\end{crlre}}
\newcommand{\exmpl}{\end{xmpl}}
\newcommand{\ermrk}{\end{rmrk}}

\newcommand{\bbn}{\mathbb{N}}

\newcommand{\bbt}{\mathbb{T}}

\newcommand{\clh}{\mathcal{H}}

\newcommand{\ovn}{\overline{\mathbb{N}}}
\newcommand{\ovz}{\overline{\mathbb{Z}}}

\title{Inverse semigroups and Sheu's groupoid for the odd dimensional quantum spheres}
\author{S.Sundar}


\begin{document}
\maketitle
\begin{abstract}
In this paper, we give a different proof of the fact that the odd dimensional
quantum spheres are groupoid $C^{*}$ algebras. We show that the $C^{*}$ algebra
$C(S_{q}^{2\ell+1})$ is generated by an inverse semigroup $T$ of partial
isometries. We show that the groupoid $\mathcal{G}_{tight}$ associated to the
inverse semigroup $T$ in \cite{Ex} is exactly the same as the groupoid
considered in \cite{Sh2}. 
\end{abstract}

\section{Introduction}
  Quantization of mathematical theories is a major theme of research today. The
theory of Quantum groups and Noncommutative
geometry are two prime examples in this program. The theory of compact quantum
groups were initiated by Woronowicz in the 
late eighties in \cite{Wor_compact_matrix}, \cite{Wor_Tannaka},
\cite{Wor_compact_quant_groups}. A main example in his theory is the 
quantum group $SU_{q}(n)$ and its homogeneous spaces. One of the problems in
noncommutative geometry is to understand how these groups fit under
Connes' formulation of NCG. Thus it becomes necessary to understand the $C^{*}$
algebra of these quantum groups.
  
 Vaksman and Soibelman studied the irreducible representations of $C^{*}$
algebra of the quantum group $SU_{q}(n)$ in \cite{Soibelman}.
Exploiting their work , Sheu in \cite{Sh1} used the theory of groupoids and
obtained certain composition sequences which is useful in 
understanding the structure of the $C^{*}$ algebra of $SU_{q}(n)$. Then in
\cite{Sh2} the question of whether $C^{*}$ algebras of these quantum homogeneous
spaces are in fact groupoid $C^{*}$ algebras was raised. In \cite{Sh2}, an affirmative answer is given for 
the quantum homogeneous space $S_{q}^{2n-1}:=SU_{q}(n)/SU_{q}(n-1)$ called the odd dimensional quantum spheres. The purpose of this paper 
is to give an alternative proof of the same result. We use the theory of inverse semigroups developed in \cite{Ex} to reconstruct the 
groupoid given in \cite{Sh2}. We believe that the proof is constructive as the groupoid in \cite{Sh2} is reconstructed 
from a combinatorial data naturally associated to $S_{q}^{2n-1}$.

Now we indicate the organisation of this paper. In the next section, we recall the basics of inverse semigroups and the groupoid 
associated to it without proofs. We refer to \cite{Ex} for proofs. In section 3, we recall the definition of the $C^{*}$ algebra
of the odd dimensional sphere $S_{q}^{2\ell+1}$ and associate a natural inverse semigroup to it. In sections 4-6, we work out the 
groupoid associated to the inverse semigroup and show that the groupoid is infact naturally isomorphic to Sheu's groupoid
constructed in \cite{Sh2}. We end the paper by showing that the reduced $C^{*}$ algebra associated to the groupoid is infact 
isomorphic to the $C^{*}$ algebra of the odd dimensional quantum sphere $S_{q}^{2\ell+1}$.

\section{ Inverse semigroups and their groupoids }
      In this section, we briefly recall the construction of the groupoid
associated to an inverse semigroup. We refer to \cite{Ex} for proofs and
details.
\begin{dfn}
  An inverse semigroup $T$ is an associative semigroup for every $s \in T$, there exists a unique element denoted $s^{*}$ such that $s^{*}ss^{*}=s^{*}$ and $ss^{*}s=s$. Then $*$ is an involution which is antimultiplicative. We say that an inverse semigroup has $0$ if there exists an element $0$ such that $0.s=s.0 =
0$ for every $s \in T$.
\end{dfn}

\subsection{The unit space of the groupoid}
 Let $T$ be an inverse semigroup with $0$. We denote the set of projections in
$T$ by $E$ i.e. $E:=\{e \in T: e=e^{*}=e^{2}\}$. Then $E$ is a commutative
semigroup containing $0$. Consider the set $\{0,1\}$ as a multiplicative
semigroup.

\begin{dfn}
 Let $T$ be an inverse semigroup with $0$ and let $E$ be its set of projections.
A character of $E$ is a nonzero map $x:E\to \{0,1\}$ such that 
\begin{enumerate}
 \item The map $x$ is a semigroup homomorphism, and 
 \item $x(0)=0$.
\end{enumerate}
We denote the set of characters of $E$ by $\hat{E}_{0}$. The set of characters
$\hat{E}_{0}$ is a locally compact Hausdorff topological space where the
topology on $\hat{E}_{0}$ is the subspace topology inherited from the topology
of $\{0,1\}^E$.
\end{dfn}

The set of characters can alternatively be described in terms of filters by
considering its support. For a character $x$, let $A_{x}:=\{e \in E: x(e)=1 \}$.
Then $A_x$ is nonempty and has the following properties
\begin{enumerate}
 \item $0 \notin A_x$,
 \item If $e \in A_x$ and $f \geq e$ then $f \in A_x$, and
 \item If $e,f \in A_x$ then $ef \in A_x$.
\end{enumerate}

A nonempty ssubset $A$ of $E$ having properties $(1),(2)$ and $(3)$ is called a
filter. Moreover, if $A$ is a filter then the indicator function $1_{A}$ is a
character. Thus there is a bijection between filters and characters. A filter is
called an ultra filter if it is maximal. By Zorn's lemma ultra filters exist.
Define
\begin{displaymath}
 \hat{E}_{\infty}:=\{x \in \hat{E}_{0}: A_x \textrm{~ is an ultrafilter } \}
\end{displaymath}
and denote its closure by $\hat{E}_{tight}$.

\subsection{ The partial action of $T$ on $\hat{E}_{0}$ }
 The inverse semigroup $T$ acts naturally on $\hat{E}_{0}$ as partial
homoemorphisms which we now explain. We let $T$ act on $\hat{E}_{0}$ on the
right as follows. 

For $x \in \hat{E}_{0}$ qnd  $s \in T$, define $(x.s)(e)=x(ses^*)$. Then
\begin{enumerate}
 \item The map $x.s$ is a semigroup homomorphism, and 
 \item $(x.s)(0)=0$.
\end{enumerate}
But $x.s$ is nonzeo if and only if $x(ss^*)=1$. For $s \in T$, define the domain
and range of $s$ as 
\begin{align*}
 D_s:&=\{x \in \hat{E}_{0}: x(ss^*)=1\} \\
 R_s:&=\{x \in \hat{E}_{0}: x(s^*s)=1\} 
\end{align*}
Note that both $D_s$ and $R_s$ are compact and open. Moreover $s$ defines a
homoemorphism from $D_s$ to $R_s$ with $s^*$ as its inverse. Also observe that
$\hat{E}_{tight}$ is invariant under the action of $T$.

\subsection{The groupoid $\mathcal{G}_{tight}$}
    Consider the transformation groupoid $\Sigma:=\{(x,s):x\in D_s\}$ with the
composition and the inversion given by:
 \begin{align*}
  (x,s)(y,t):&=(x,st) \textrm{~if~} y=x.s\\
   (x,s)^{-1}:&=(x.s,s^*) 
 \end{align*}
Define an equivalence relation $\sim$ on $\Sigma$ as $(x,s)\sim (y,t)$ if $x=y$
and if there exists an $e \in E$ such that $x \in D_e$ for which $es=et$. Let
$\mathcal{G}=\Sigma/\sim$. Then $\mathcal{G}$ is a groupoid as the product and
the inversion respects the equivalence relation $\sim$. Now we describe a
toplogy on $\mathcal{G}$ which makes $\mathcal{G}$ into a topological groupoid.

For $s \in T$ and $U$ an open subset of $D_s$, let $\theta(s,U):=\{[x,s]: x \in
U \}$. We refer to \cite{Ex} for the proof of the following two propositions. We
denote $\theta(s,D_s)$ by $\theta_s$. Then $\theta_s$ is homeomorphic to $D_s$
and hence is compact, open and Hausdorff.

\begin{ppsn}
 The collection $\{\theta(s,U): s \in T, U \textrm{~open in ~} D_s \}$ forms a
basis for a topologyon $\mathcal{G}$. The groupoid $\mathcal{G}$ with this
topology is a topological groupoid whose unit space can be identified with
$\hat{E}_{0}$.
Also one has the following.
\begin{enumerate}
 \item For $s,t \in T$, $\theta_s\theta_t=\theta_{st}$,
 \item For $s \in T$, $\theta_{s}^{-1}=\theta_{s^{*}}$, and
 \item The set $\{1_{\theta_s}: s \in T\}$ generates the $C^*$ algebra
$C^*(\mathcal{G})$.
\end{enumerate}
\end{ppsn}

We define the groupoid $\mathcal{G}_{tight}$ to be the reduction of the groupoid
$\mathcal{G}$ to $\hat{E}_{tight}$.

\section{The odd dimensional quantum spheres}
 Before we proceed let us fix some notations. Throughout we assume that $q \in
(0,1)$ and $\ell$ is a positive integer. We denote the set of non-negative integers by $\mathbb{N}$. We let $\clh_{\ell}$ to denote the
Hilbert space $\ell^{2}(\mathbb{N})^{\otimes \ell} \otimes
\ell^{2}(\mathbb{Z})$. We denote the left shift operator on
$\ell^{2}(\mathbb{N})$ by $S$ and the right shift on $\ell^{2}(\mathbb{Z})$ by
$t$. The number operator on $\ell^{2}(\mathbb{N})$ is denoted by $N$.

In this section we recall a few well known facts about the $C^*$ algebra of the
odd dimensional quantum spheres.
The $C^*$-algebra $C(S_q^{2\ell+1})$ of the quantum
sphere $S_q^{2\ell+1}$
is the universal $C^*$-algebra generated by
elements
$z_1, z_2,\ldots, z_{\ell+1}$
satisfying the following relations (see~\cite{Hong-Sym1}):
\bean
z_i z_j & =& qz_j z_i,\qquad 1\leq j<i\leq \ell+1,\\
z_i^* z_j & =& q z_j z_i^* ,\qquad 1\leq i\neq j\leq \ell+1,\\
z_i z_i^* - z_i^* z_i +
(1-q^{2})\sum_{k>i} z_k z_k^* &=& 0,\qquad \hspace{2em}1\leq i\leq \ell+1,\\
\sum_{i=1}^{\ell+1} z_i z_i^* &=& 1.
\eean
 Note that for $\ell=0$, the $C^*$-algebra
$C(S_q^{2\ell+1})$ is the algebra of continuous functions
$C(\bbt)$ on the torus and for $\ell=1$, it is $C(SU_q(2))$.

Let $Y_{k,q}$ be the following operators on $\clh_\ell$:
\be\label{eq:ykq}
 Y_{k,q}=\begin{cases}
 \underbrace{q^N\otimes\ldots\otimes q^N}_{k-1 \mbox{
copies}}\otimes
      \sqrt{1-q^{2N}}S^*\otimes 
   \underbrace{I \otimes\cdots\otimes I}_{\ell+1-k \mbox{ copies}}, & \mbox{ if
} 1\leq k\leq \ell,\cr
   &\cr
    \underbrace{q^N\otimes\cdots\otimes q^N}_{\ell \mbox{ copies}}
       \otimes t, &  \mbox{ if } k=\ell+1.
         \end{cases}
\ee
Then $\pi_\ell:z_k\mapsto Y_{k,q}$ gives a faithful representation
of $C(S_q^{2\ell+1})$ on $\clh_\ell$ for $q\in(0,1)$ 
(see lemma~4.1 and remark~4.5, \cite{Hong-Sym1}). We let $Y_{k,0}$ to denote the
limit of the operators $Y_{k,q}$ as $q$ tends to zero. The formula for
$Y_{k,0}$ is again the same as that of $Y_{k,q}$ where $q^N$ stands for the rank
one projection $p=|e_0\rangle\langle
e_0|$. 

Consider the unitary operator $U$  on $\clh$ defined by
$U(e_{m,z})=e(m,z+\sum_{i=1}m_{i})$. Define $Z_{k,q}:=UY_{k,q}U^{*}$ for $q \in
[0,1]$.
The representation $z_{k} \to Z_{k,q}$ of $C(S_{q}^{2\ell+1})$ is the one
considered in \cite{Sh1} and in \cite{Sh2}.
Let $A_{\ell}(q)$ be the image of $C(S_{q}^{2\ell+1})$ under this representation
i.e. $A_{\ell}(q)$ is the $C^{*}$ algebra 
generated by $Z_{k,q}$. We refer to \cite{SunPal} for the proof of the following
proposition.

\begin{ppsn}
For any $q\in(0,1)$, one has $A_\ell(0)=A_\ell(q)$.
\end{ppsn}

From now on, we simply denote $Z_{k,0}$ by $Z_{k}$. Note that $Z_{k}$'s are in
fact partial isometries.
Let us introduce more notations. For $m,n \in \mathbb{N}^{\ell}$ and $r\in
\mathbb{Z}$, Let $B_{k}(r,m,n)$ be defined
as follows:
     \be\label{eq:Bk}
 B_{k}(r,m,n)=\begin{cases}
 \underbrace{S^{*m_{1}}pS^{n_{1}}\otimes S^{*m_{2}}pS^{n_{2}}\otimes
\ldots\otimes S^{*m_{k-1}}pS^{n_{k-1}}}_{k-1 \mbox{
copies}}\otimes
      S^{*m_{k}}S^{n_{k}}\otimes 1 \otimes t^{(\sum_{i=1}^{k}(m_{i}-n_{i}))}
    &  \mbox{ if
} 1\leq k\leq \ell,\cr
   &\cr
    \underbrace{S^{*m_{1}}pS^{n_{1}}\otimes S^{*m_{2}}pS^{n_{2}}\otimes
\ldots\otimes S^{*m_{\ell}}pS^{n_{\ell}}} \otimes
t^{r+\sum_{i=1}^{\ell}(m_{i}-n_{i})}, &  \mbox{ if } k=\ell+1.
         \end{cases}
\ee

Note the following commutation relations.

If $i<j$ then \[
 B_{i}(r,m,n)B_{j}(r^{'},m^{'},n^{'})=
\delta_{(n_{1},n_{2},\cdots,n_{i-1}),(m^{'}_{1},m_{2}^{'},\cdots
m_{i-1}^{'})}1_{[n_{i},\infty)}(m_{i}^{'})B_{k}(r^{'},m^{''},n^{'})            
              \]
 where
$m^{''}=(m_{1},m_{2},\cdots,m_{i-1},m_{i}+m_{i}^{'}-n_{i},m^{'}_{i+1},\cdots,m^{
'}_{\ell})$.

If $i\leq \ell$ and $n_{i} \leq m_{i}^{'}$ then \[
B_{i}(r,m,n)B_{i}(r^{'},m^{'},n^{'})= 
\delta_{(n_{1},n_{2},\cdots,n_{i-1}),(m_{1}^{'},m_{2}^{'},\cdots,m_{i-1}^{'})}B_
{i}(r^{'},m^{''},n^{'})
                                                \]
 where
$m^{''}:=(m_{1},m_{2},\cdots,m_{i-1},m_{i}+m_{i}^{'}-n_{i},m_{i+1}^{'},\cdots
m_{l}^{'})$.

If $i\leq \ell$ and $n_{i} < m_{i}^{'}$ then \[
                                        B_{i}(r,m,n)B_{i}(r^{'},m^{'},n^{'})= 
\delta_{(n_{1},n_{2},\cdots,n_{i-1}),(m_{1}^{'},m_{2}^{'},\cdots,m_{i-1}^{'})}B_
{i}(r^{'},m,n^{''})\]

where
$n^{''}:=(n_{1}^{'},n_{2}^{'},\cdots,n_{i-1}^{'},n_{i}^{'}+n_{i}-m_{i}^{'},n_{
i+1}^{'},\cdots,n_{\ell}^{'})$.

Finally,
$B_{\ell+1}(r,m,n)B_{\ell+1}(r^{'},m^{'},n^{'})=\delta_{n,m^{'}}B_{\ell+1}(r+r^{
'},m,n^{'})$. Also $B_{i}(r,m,n)^{*}:=B_{i}(-r,n,m)$. It is clear from the above
commutation relations that 
the set $T:=\{0\}\bigcup \{B_{i}(r,m,n): 1\leq i \ leq \ell+1, r \in
\mathbb{Z},m,n \in \mathbb{N}^{k}\}$ is an inverse semigroup of partial
isometries. 

\begin{ppsn}
\label{Inverse semigroup}
 The set $T:=\{0\}\bigcup \{B_{i}(r,m,n): 1\leq i \leq \ell+1, r \in
\mathbb{Z},m,n
\in \mathbb{N}^{k}\}$ is an inverse semigroup of partial isometries. Moreover
$T$ is generated by $\{Z_{i}:1\leq i \leq \ell+1\}$. 
\end{ppsn}
\textit{Proof.} As already observed $T$ is an inverse semigroup of partial
isometries. Let $e_{i}$ be the $\ell$ tuple which is $1$ on the $i$th coordinate
and zero elsewhere.
Then $Z_{k}:=B(0,e_{k},0)$ for $k \leq \ell$ and $Z_{\ell+1}=B_{\ell+1}(1,0,0)$.
Thus $Z_{k}$'s are in $T$. Moreover,
 \begin{align*}
        0 &=Z_{1}^{*}Z_{2} \\
B_{i}(r,m,n)&=Z_{1}^{m_{1}}Z_{2}^{m_{2}}\cdots
Z_{i}^{m_{i}}Z_{i}^{*n_{i}}Z_{i-1}^{*n_{i-1}}\cdots Z_{1}^{*n_{1}}\mbox{~if~}
i\leq \ell \\
B_{\ell+1}(r,m,n)&=Z_{1}^{m_{1}}Z_{2}^{m_{2}}\cdots
Z_{\ell}^{m_{\ell}}(Z_{\ell+1}^{r_{+}}Z_{\ell+1}^{r_{-}})Z_{\ell}^{*n_{\ell}}Z_{
\ell-1}^{*n_{\ell-1}}\cdots Z_{1}^{*n_{1}}.
\end{align*}
where $r_{+}$ and $r_{-}$ denote the positive and negative parts of $r$. Thus
every element in $T$ is a word in $Z_{i}$'s.
This completes the proof. \hfill $\Box$

\section{The tight characters for the inverse semigroup $T$}
In this section, we describe the tight characters of the inverse semigroup $T$
defined in Proposition \ref{Inverse semigroup}.
The set of projections of $T$ is denoted by $E$ and the set of characters of $E$
by $\hat{E}_{0}$. Consider the one point compactification
$\overline{\mathbb{N}}:=\mathbb{N}\cup \{\infty\}$ of $\mathbb{N}$. We denote
the projection from $\ovn^{\ell}$ onto the first $r$ components by $\pi_r$.

Let $p_{i}(m):=B_{i}(0,m,m)$. Then $E:=\{0\}\cup \{p_{i}(m):1 \leq i \leq
\ell+1, m \in \mathbb{N}^{k}\}$. First observe that if $\Lambda$ is a
subsemigroup of $E$ not containing $0$
then the set \[
              A_{\Lambda}:=\{f \in E: f \geq e \mbox{~for some ~} e \in \Lambda
\}
             \]
 is a filter.

Let $k \in \ovn^{\ell}$ be given. Let $r$ be the least positive integer for
which $k_{r+1}=\infty$. ( For an $\ell$ tuple $k$, we set $k_{\ell+1}=\infty$).
Define $\Lambda_k=\{p_{r+1}(\pi_r(k),n): n \in \bbn^{\ell} \}$. It is easy to
see that $\Lambda_k$ is a subsemigroup of $E$ not containing $0$. Then
$A_{\Lambda_k}$ is a filter and thus gives rise to a character. We denote the
character associated to $A_{\Lambda_k}$ by $\phi(k)$. The following lemma gives
a closed formula for $\phi(k)$.

\begin{lmma}
\label{formula for phi}
 Let $k \in \ovn^{\ell}$ be given. The character $\phi(k)$ is given by 
\begin{displaymath}
\begin{array}{lll}
\phi(k)(p_i(m))&:=\left\{\begin{array}{lll}
                                
\delta_{\pi_{i-1}(m),\pi_{i-1}(k)}1_{[0,k_i]}(m_i) & if & i \leq \ell \\
                                 \delta_{m,k} & if & i=\ell+1 
                                 \end{array} \right. 
\end{array}
\end{displaymath}
\end{lmma}
\textit{Proof.} Let $r$ be the least positive integer for which
$k_{r+1}=\infty$. Observe  $p_{i}(m) \geq p_{r+1}(\pi_r(k),n)$ for some $n$ if
and only if $i \leq r+1$, $\pi_{i-1}(m)=\pi_{i-1}(k)$ and $m_{i} \leq k_{i}$.
Now the proof follows. \hfill $\Box$

An immediate consequence of the above lemma is that the map $\phi:\ovn^{\ell}\to
\hat{E}_0$ is continuous. In the next proposition we show that the image of
$\phi$ is exactly $\hat{E}_{\infty}$.

\begin{ppsn}
\label{tight characters}
 The image of the map $\phi$ is infact $\hat{E}_{\infty}$. 
\end{ppsn}
\textit{Proof.} First we show that the image of $\phi$ is contained in
$\hat{E}_{\infty}$. Let $k \in \ovn^{\ell}$ be given and let $r$ be the least
nonnegative integer for which $k_{r+1}=\infty$.
Recall that 
\[
 \Lambda_{k}:=\{p_{r+1}(\pi_{r}(k),n): n \in \bbn^{\ell-r}\}
\]
We denote $\phi(k)$ by $x$. We claim that $A_{x}$ is an ultrafilter. Suppose
that there exists a character say $y$ such that $A_{x} \subset A_{y}$.
Then we need to show that $x=y$ or $A_{x}=A_{y}$. Since $x=1$ on $\Lambda_{k}$
and $A_{x}\subset A_{y}$, it follows that $y=1$ on $\Lambda_{k}$.
If $\pi_{r}(m)\neq \pi_{r}(k)$ then $p_{r+1}(\pi_{r}(m),v)$ is orthogonal to
$\Lambda_{k}$. Hence $x$ and $y$ vanishes on $p_{r+1}(\pi_{r}(m),v)$.
Thus $x(p_{r+1}(m))=y(p_{r+1}(m))$ for every $m \in \mathbb{N}^{\ell}$..

Now let $i>r+1$ be given. Let $m \in \mathbb{N}^{\ell}$. Choose an $\ell$ tuple
$n$ such that $\pi_{r}(n)y=\pi_{r}(k)$ and $n_{r+1}>m_{r+1}$.
Then $p_{i}(m)$ and $p_{r+1}(n)$ are orthogonal. But $x=y=1$ on $p_{r+1}(n)$.
Thus $x$ and $y$ vanishes on $p_{i}(m)$.

Now let $i\leq r$ and $m \in \mathbb{N}^{\ell}$. If $m_{i}>k_{i}$ then
$p_{i}(m)$ is orthogonal to $\Lambda_{k}$. But $x=y=1$ on $\Lambda_{k}$.
Thus $x$ and $y$ vanishes on $p_{i}(m)$ if $m_{i}>k_{i}$. Now suppose $m_{i}\leq
k_{i}$. If $\pi_{i-1}(m) \neq \pi_{i-1}(k)$ then 
$p_{i}(m)$ is again orthogonal to $\Lambda_{k}$ and thus $x$ and $y$ vanishes on
$p_{i}(m)$. Consider now the case where
$m_{i} \leq k_{i}$ and $\pi_{i-1}(m)=\pi_{i-1}(k)$. Then $x(p_{i}(m))=1$ by
definition and since $A_{x} \subset A_{y}$, it follows that
$y(p_{i}(m))=1$. Thus we have shown that $x(p_{i}(m))=y(p_{i}(m))$ for every $i$
and $m$. Hence $x=y$ or in otherwords
$A_{x}$ is an ultrafilter. This proves that $\phi(\ovn^{\ell})$ is contained in
$\hat{E}_{\infty}$.

Now let us prove that $\hat{E}_{\infty}$ is contained in the range of $\phi$.
Let $x\in \hat{E}_{\infty}$ be given. Let $r$ be the largest nonnegative integer
for which there
exists a $k^{'}$ such that $x=1$ on $p_{r+1}(k^{'})$. Choose $k$ such that
$\pi_{r}(k)=\pi_{r}(k^{'})$ and $k_{r+1}=\infty$
We claim that $A_{x} \subset A_{\phi(k)}$. Then the maximality of $A_{x}$ forces
the equality $x=\phi(k)$.

Let $i \leq r+1$ be given. Consider an $\ell$ tuple $m$ such that either
$m_{i}>k_{i}$ or $\pi_{i-1}(m)\neq \pi_{i-1}(k)$.
Then $p_{i}(m)$ is orthogonal to $p_{r+1}(k^{'})$. Hence $x(p_{i}(m))=0$ if
either $m_{i}>k_{i}$ or $\pi_{i-1}(m) \neq \pi_{i-1}(k)$.
Also $x$ vanishes on $p_{i}(m)$ if $i>r+1$ by the choice of $r$. Thus we have
shown $A_{\phi(k)}^{c} \subset A_{x}^{c}$. Hence
$A_{x}$ is contained in $A_{\phi(k)}$. Since $A_{x}$ is maximal, it follows that
$x=\phi(k)$. This completes the proof.
\hfill $\Box$

\begin{crlre}
 The set $\hat{E}_{\infty}$ is compact and  $\hat{E}_{tight}=\hat{E}_{\infty}$.
\end{crlre}
\textit{Proof.} The proof follows from the fact that $\phi$ is continuous,
$\ovn^{\ell}$ is compact and from Proposition \ref{tight characters}.
\hfill $\Box$.

Now define an equivalence relation on $\ovn^{\ell}$ as follows: 
\[
 k \sim k^{'} \text{~ if there exists } r\geq 0 \text{~such that }
\pi_{r}(k)=\pi_{r}(k^{'}) \text{~and~} k_{r+1}=k^{'}_{r+1}.
\]
We show that $\hat{E}_{tight}$ is homeomorphic to the qutient space
$\ovn^{\ell}/\sim$ in the next proposition.
\begin{ppsn}
\label{isomorphism of tight characters}
 The map $\phi:\ovn^{\ell} \to \hat{E}_{tight}$ factors through the quotient
$\ovn^{\ell}/\sim$ to give a map $\tilde{\phi}:\ovn^{\ell}/\sim \to
\hat{E}_{tight}$. Also the maps
$\tilde{\phi}$ is a homeomorphism.
\end{ppsn}
\textit{Proof.} It is clear from the definition and from Lemma \ref{formula for
phi} that $\phi$ factors through the 
quotient to give a map $\tilde{\phi}$. Since $\phi$ is continuous, it follows
that $\tilde{\phi}$ is continuous. We now show that
$\tilde{\phi}$ is one-one.

Let $k,k^{'}$ be such that $\phi(k)=\phi(k^{'})$. Let $r_{k}$ (resp.
$r_{k^{'}}$) be the least nonnegative integer
for which $k_{r_{k}+1}=\infty$ (resp. $k^{'}_{r_{k^{'}}+1}=\infty$). Then
$r_{k}$ is the last integer for which there exists an 
$m$ such that $\phi(k)$ is $1$ on $p_{r_{k}+1}(m)$. Thus $r_{k}=r_{k^{'}}$.
Moreover $\phi(k)$ is $1$ on $p_{r_{k}+1}(\pi_{r_{k}}(k),u)$.
Thus $\phi_{k^{'}}$ is $1$ on $p_{r_{k+1}}(\pi_{r_{k}}(k),u)$. Now Lemma
\ref{formula for phi} implies that $\pi_{r_{k}}(k)=\pi_{r_{k}}(k^{'})$.
Hence $k \sim k^{'}$. This proves that $\tilde{\phi}$ is one-one.

Now Proposition \ref{tight characters} implies that $\tilde{\phi}$ is onto. As
$\ovn^{\ell}/\sim$ is compact and $\hat{E}_{tight}$ is Hausdorff,
it follows that $\tilde{\phi}$ is infact a homeomorphism. This completes the
proof. \hfill $\Box$

\section{Sheu's groupoid}
 In this section, we recall the groupoid for the odd dimensional quantum spheres
$S_{q}^{2\ell+1}$ described in \cite{Sh2}.
Consider the transformation groupoid $\mathbb{Z} \times (\mathbb{Z}^{\ell}
\times \ovz^{\ell})$ where $\mathbb{Z}^{\ell}$ acts on 
$\ovz^{\ell}$ by translation. Let $\mathcal{F}$ be the restriction of the
transformation groupoid to $\ovn^{\ell}$. Define
\[
 \Sigma:=\{(z,x,w) \in \mathcal{F}: w_{i}=\infty \Rightarrow
x_{i+1}=\cdots=x_{\ell}=0 \text{~and } z= -\sum_{j=1}^{i}x_{j} \}
\]
Then $\Sigma$ is an open subgroupoid of $\mathcal{F}$. Define an equivalence
relation $\sim$ on $\Sigma$ as follows:
\[
 (z,x,w_{1},w_{2},\cdots,w_{i-1},\infty,\cdots) \sim
(z,x,w_{1},w_{2},\cdots,w_{i-1},\infty,\infty,\cdots,\infty).
\]
Let $\mathcal{G}:=\Sigma/\sim$. Then the multiplication and the inversion on
$\Sigma$ factors through the equivalence relation
making $\mathcal{G}$ into a groupoid. When $\mathcal{G}$ is given the quotient
topology, it becomes a topological groupoid which is the groupoid described in
\cite{Sh2}.

\section{The groupoid $\mathcal{G}_{tight}$ of the inverse semigroup $T$}
  In this section, we show that the groupoid $\mathcal{G}_{tight}$ of the
inverse semigroup $T$ is isomorphic
with Sheu's groupoid described in the previous section. For an $\ell$ tuple $m$,
we set $m_{\ell+1}=\infty$. We define
a map $\psi:\Sigma \to \mathcal{G}_{tight}$ as follows:

Let $(z,x,w) \in \Sigma$ be given. Let $r$ be the least nonnegative integer for
which $w_{r+1}=\infty$. Then $\psi$ on
$(z,x,w)$ is given by 
\[
   \psi(z,x,w):=[(\phi(w), B_{r+1}(t,m,n)]
\]
where $t,m,n$ are given by $t:=z+\sum_{j=1}^{r}x_{j}$,
$m:=(w_{1},w_{2},\cdots,w_{r},|x_{r+1}|,0,\cdots,0)$ and
 $n:=(x_{1}+w_{1},x_{2}+w_{2},\cdots,x_{r}+w_{r},x_{r+1}+|x_{r+1}|,0,\cdots,0)$.
Observe that $\psi$ is well defined 
as $\pi_{r}(w)=\pi_{r}(m)$ and $w_{r+1}=\infty$. 

Let us introduce the following notation. For $m,n \in \mathbb{N}^{\ell}$ and
$0\leq r \leq \ell$, let  $A_{r+1}(m,n):=S^{*m_{1}}pS^{n_{1}}\otimes \cdots
\otimes S^{*m_{r}}pS^{n_{r}}\otimes  S^{*m_{r+1}}S^{n_{r+1}} \otimes 1$. We
consider $A_{r+1}(m,n)$ as an operator on $\ell^{2}(\mathbb{N}^{\ell})$.

\begin{ppsn}
 The map $\psi$ is continuous and $\psi$ factors through the equivalence
relation $\sim$. Let $\tilde{\psi}$ be the induced map
from $\mathcal{G} \to \mathcal{G}_{tight}$. Then $\tilde{\psi}$ is infact a
groupoid isomorphism.
\end{ppsn}
\textit{Proof.} First we show that $\psi$ factors through the equivalence
relation. Let $(z,x,w)\sim (z,x,w^{'})$ and let $r$ (resp. $r^{'}$) 
be the least nonnegative integer for which $w_{r+1}=\infty$ (resp
$w^{'}_{r^{'}+1}=\infty$). By definition $r=r^{'}$ and
$\pi_{r}(w)=\pi_{r}(w^{'})$.
Then by Proposition \ref{isomorphism of tight characters}, it follows that
$\phi(w)=\phi(w^{'})$. Since the definition of $\psi$ involves only the first
$r$ components of $w$,  $\psi(z,x,w)=\psi(z,x,w^{'})$. This proves that
$\tilde{\psi}$ is well defined.

\underline{The map $\tilde{\psi}$ is one-one.} Suppose that
$\psi(z,x,w)=\psi(z^{'},x^{'},w^{'})$. Again let $r$ and $r^{'}$ be 
the least nonnegative integer for which $w_{r+1}$ and $w_{r^{'}+1}$ are both
$\infty$. Then $r$ is the largest integer for which there exists an $m$ such
that 
$\phi(w)$ is $1$ on $p_{r+1}(m)$. Since $\phi(w)=\phi(w^{'})$, we have $r=r^{'}$
and $\pi_{r}(w)=\pi_{r}(w^{'})$. As $\psi(z,x,w)=\psi(z^{'},x^{'},w^{'})$, it
follows that
there exists a projecti on $e$ such that $\phi(w)(e)=1$ and $e(t^{z}\otimes
A_{r+1}(m,n))=e(t^{z^{'}} \otimes A_{r+1}(m^{'},n^{'})$.
But $\phi(w)(e)=1$ implies that $e \geq p_{r+1}(\pi_{r}(w),0)$. Hence we can
choose $e$ to be $p_{r+1}(\pi_{r}(w),0)$.
Thus it follows that $z=z^{'}$ and $A_{r+1}(m,n)=A_{r+1}(m^{'},n^{'})$. Thus
$m=m^{'}$ and $n=n^{'}$ which in turn implies $x_{i}=x_{i}^{'}$ for $i \leq
r+1$.
Since $(z,x,w) \in \Sigma$, it follows that $x_{i}=x_{i}^{'}=0$ for $i \geq
r+2$. Thus we have shown that $(z,x,w) \sim (z^{'},x^{'},w^{'})$. 
Hence $\tilde{\psi}$ is one-one.

\underline{The map $\tilde{\psi}$ is onto.} First note that if $a-b=c-d$ then
there exists a projection $e=S^{*(b+c)}S^{b+c}$ such that
$eS^{*b}S^{a}=eS^{*d}S^{c}$. Hence in the definition of $\psi$ we can change the
$r+1$th components of $m$ and $n$ such that
$n_{r+1}-m_{r+1}=x_{r+1}$. Let $[(\phi(w),B_{i}(s,m,n))]$ be an element in
$\mathcal{G}_{tight}$. Let $r$ be the first 
nonnegative integer for which $w_{r+1}=\infty$. Then $i \leq r+1$. By
premultiplying by $p_{r+1}(\pi_{r}(w),0)$ we can assume that
$i=r+1$ and $m$ is such that $\pi_{r}(w)=\pi_{r}(m)$. Now if $r \leq \ell-1$
then ,for $z=\sum_{j}(m_{j}-n_{j})$, $x$ such
that$\pi_{r+1}(x)=\pi_{r+1}(n)-\pi_{r+1}(m)$
and $x_{i}=0$ for $i \geq r+2$, $\psi(z,x,w)=[(\phi(w),B_{r+1}(s,m,n))]$. If
$r=\ell$ then with $z=s+\sum_{j=1}^{\ell}(m_{j}-n_{j}$ and $x_{j}=n_{j}-m_{j}$
one has 
$\psi(z,x,w)=[(\phi(w),B_{\ell+1}(s,m,n))]$. This proves that $\tilde{\psi}$ is
onto.

\underline{The map $\tilde{\psi}$ is continuous.} Let $(z^{n},x^{n},w^{n})$ be a
sequence in $\Sigma$ converging to $(z,x,w) \in \Sigma$.
Let $r$ be the least nonnegative integer for which $w_{r+1}=\infty$. Then
eventually $(z^{n},x^{n},\pi_{r}(w^{n}))$ coincides with $(z,x,\pi_{r}(w))$.
Suppose that $\theta(s,U)$ is an open set containing $\psi(z,x,w)$. Without loss
of generality 
we can assume that $s:=t^{z}\otimes A_{r+1}(m,n)$ where $m,n$ are defined as in
the definition
of $\psi$. Since $U$ is an open set containing $\phi(w)$ and as $\phi$ is
continuous, it follows
that $\phi(w^{n}) \in U$ eventually. Let $r_{n}$ be the least nonnegative
integer for which $w^{n}_{r_{n}+1}=\infty$.
Then $r_{n} \geq r$. Let $m^{n}$ and $n^{n}$ be as in the definition of $\psi$
for $(z,x,w^{n})$. If $e_{n}:=p_{r_{n}+1}(m^{n})$, then $\phi(w^{n}) \in
D_{e_{n}}$ 
and $e_{n}(t^{z}\otimes A_{r_{n}+1}(m^{n},n^{n}))=e_{n}(t^{z}\otimes
A_{r+1}(m,n))$ eventually.
Thus eventually $\psi(z^{n},x^{n},w^{n})=[(\phi(w^{n}),s)] \in \theta(s,U)$.
This proves that $\tilde{\psi}$ is continuous.

We leave it to the reader to check that $\tilde{\psi}$ is a homeomorphism and it
is infact a groupoid homomorphism.
This completes the proof. \hfill $\Box$

\section{Isomorphism between $C(S_{q}^{2\ell+1})$ and
$C^{*}_{red}(\mathcal{G})$}
 We complete the discussion by showing that $C(S_{q}^{2\ell+1})$ is isomorphic
to the reduced $C^{*}$ algebra of the 
groupoid $\mathcal{G}$ where $\mathcal{G}$ is the Sheu's groupoid considered in
Section 5. We also identify $\mathcal{G}$ with $\mathcal{G}_{tight}$.
Let $r$ denote the range map and let $\mathcal{G}^{0}:=r^{-1}(\phi(0))$. Then
$\mathcal{G}^{0}:=\{(z,x,0):x_{i} \geq 0\}$.
Thus $L^{2}(\mathcal{G}^{0})$ is naturally isomorphic to
$\ell^{2}(\mathbb{N}^{\ell})\otimes \ell^{2}(\mathbb{Z})$.
Consider the representation $\pi$ of $C^{*}_{red}(\mathcal{G})$ on
$L^{2}(\mathcal{G}^{0}$ defined by: 

For $f \in C_{c}(\mathcal{G})$ and $\xi \in L^{2}(\mathcal{G}^{0})$, 
\[
 (\pi(f)\xi)(\gamma)=\sum_{\gamma_{1} \in
\mathcal{G}^{0}}f(\gamma^{-1}\gamma_{1})\xi(\gamma_{1})
\]
The representation $\pi$ is equivalent to $Ind (\delta_{0})$ given in
\cite{Muhly} and it is faithful since the largest
open set which is also invariant for which $\delta_{0}(U)=0$ is the empty set.
(Refer to \cite{Muhly} for more details.)

Let $\Sigma_{fin} \subset \Sigma$ be the finite part i.e.
$\Sigma_{fin}:=\{(z,x,w): w_{i} < \infty\}$. Then $\Sigma_{fin}$ is a subset of 
$\mathcal{G}$. For $1 \leq k \leq \ell+1$ , denote the set $\theta_{Z_{k}^{*}}$
by $\theta_{k}$. Now it is easily verifiable that
$\theta_{k} \cap
\Sigma_{fin}:=\{(-1,e_{k},(0,0,\cdots,0,w_{k},\cdots,w_{\ell})\}$ where $e_{k}$
is the $\ell$ tuple which is  $1$ 
at the $kth$ place and zero elsewhere. From this observation it is easy to show
that $\pi(1_{\theta_{k}})=Z_{k}^{*}$.
Since $1_{\theta_{k}}$ generate the $C^{*}$ algebra
$C_{red}^{*}(\mathcal{G}_{tight})$, it follows that $C(S_{q}^{2\ell+1})$ is
isomorphic to $C_{red}^{*}(\mathcal{G})$. 
\bibliography{references}
\bibliographystyle{plain}
\noindent
\begin{footnotesize}\textbf{Acknowledgement}: 
This work was undertaken when I was a postdoc at Universite de Caen, France from Sep 2010- Jan 2011. I thank Prof. Emmanuel Germain for his support and encouragement during my stay at Caen. I would also like to acknowledge CNRS for the financial support.\\
\end{footnotesize}
\noindent{\sc email:} \texttt{sundarsobers@gmail.com}\\
         {\footnotesize Visiting Institute of Mathematical Sciences, Chennai.}
 
\end{document}